\documentclass[9pt,twocolumn]{article}
\usepackage{setspace}
\usepackage{amsmath, amssymb,paralist,amsthm, bbm}
\allowdisplaybreaks

\def\1{\mathbbm{1}}

\def\PC{\mathcal{P}}

\def\FC{\mathcal{F}}

\def\S{\mathbf{S}}
\def\C{\mathbf{C}}
\def\E{\mathbf{E}}

\def\P{\mathbf{P}}
\def\R{\mathbf{R}}

\def\1{\mathbf{1}}

\def\al{\alpha}
\def\be{\beta}

\newcommand{\si}{\sigma}

\newcommand{\De}{\Delta}

\def\B{\mathbf{B}}
\def\D{\mathbf{D}}
\def\E{\mathbb{E}}
\def\R{\mathbf{R}}

\def\M{\mathcal{M}}
\def\O{\mathcal{O}}
\def\P{\mathcal{P}}
\def\U{\mathcal{U}}

\newtheorem{prop}{Proposition}[section]
\newtheorem{theorem}{Theorem}[section]

\newtheorem{example}{Example}[section]
\newtheorem{remark}{Remark}[section]

\numberwithin{equation}{section}

\begin{document}

\centerline{\Large \bf Sensitivity analysis for HJB }
\smallskip
\centerline{\Large \bf equations with an application to}
\smallskip
\centerline{\Large \bf coupled  backward-forward systems}
\bigskip
\centerline{\bf Vassili Kolokoltsov\footnote{Department of Statistics, University of Warwick, Coventry, CV4 7AL, UK, v.kolokoltsov@warwick.ac.uk},  Wei Yang\footnote{Department of Mathematics and Statistics, University of Strathclyde Glasgow, G1 1XH, UK, w.yang@strath.ac.uk}}
\bigskip

\begin{abstract}
In this paper, we analyse Lipschitz continuous dependence of the solution to Hamilton-Jacobi-Bellman equations on a functional parameter. This sensitivity analysis not only has the interest on its own, but also is important for the mean field games methodology, namely for solving a coupled system of backward-forward equations. We show that the unique solution to a Hamilton-Jacobi-Bellman equation and its spacial gradient are Lipschitz continuous uniformly with respect to the functional parameter.  In particular, we provide verifiable criteria for the so-called feedback regularity condition. Finally as an application, we show how the sensitive results are used to solved the coupled system of backward-forward equations.
\end{abstract}


\section{Introduction and motivation}

Partial differential equations are used to model multidimensional dynamical systems and to describe a wide variety of phenomena such as sound, heat,  fluid flow and quantum mechanics. 
Sensitivity analysis for multidimensional dynamical systems governed by partial differential equations has been a growing interest in recent years,  since sensitivity results give us certain understanding of the stability of dynamical systems and prediction of their phase transitions if any.  
So far sensitivity results have had a wide-range of applications in science and engineering, including optimization, parameter estimation, model simplification, optimal control and experimental design. 

Recent progress in sensitivity analysis can be found, e.g. in \cite{S2005} for Burger's equation, \cite{RT2003} for Navier-Stokes equation, \cite{AGMR2010, M2002, M2011,MT2000} for elliptic and parabolic equations, \cite{ Bai,Ko06a,  ManNorrisBai} for nonlinear kinetic equations, and references therein. This work at hand contributes to the presently ongoing investigation of sensitivity analysis for partial differential equations and focuses on Hamilton-Jacobi-Bellman (HJB) equations. 

The Hamilton-Jacobi-Bellman equation is a central object of the stochastic control theory.  
The solution to the HJB equation is the so-called {\it value function} which gives the optimal payoff for a given stochastic dynamical system associated with certain payoff function. 
In the context of pension savings management, Macov\'a and Sevcovic \cite{MS2010} studied sensitivity analysis of the solution to a HJB equation with respect to real parameters such as the percentage of salary transferred to a pension fund and the saver's risk aversion. 
In this paper, we study the sensitivity of the solution to HJB equations with respect to a {\it functional} parameter. 

This work is motivated by the study of the mean field games methodology, which is a recently developed research area and has attracted massive attention. The mean field games methodology was developed independently by J.-M. Lasry and P.-L. Lions, see \cite{LL2006, GLL2010} and
video lectures \cite{L}, and by M. Huang, R.P. Malham\'e and P. Caines, see \cite{HCM3, HCM07, HCM10,Hu10}.  
The Mean field games have exhibit their strong power in solving complex dynamical systems and 
have been applied  in many areas, to name a few, including growth theory in economics \cite{LLG2010},  limit order books modelling in finance \cite{LLLL2014}, environmental policy design \cite{LST2010}.   See \cite{GS2014}  for a survey on mean field games models.

The very core idea of the mean field games methodology is the characterisation and approximation of interacting N-player stochastic games by two coupled forward-backward equations, with one of these two equations being a HJB equation.
Solving this system of two coupled equations requires a so-called {\it feedback regularity property} of the feedback control associated to the HJB equation, see e.g. the condition $(37)$ in proving Theorem 10 of \cite{HCM3}.  The feedback regularity property of a control strategy is that the control strategy has a  Lipschitz continuous dependence on certain observable parameter, see the condition \eqref{FBR} in the following for its' mathematical description in the context of mean field games.
However the fact, that in \cite{HCM3}  this critical feedback regularity property was directly assumed to hold in the model without stating verifiable conditions for the feedback regularity property,  makes it impossible to apply this powerful tool in solving real life problems. 
In order to facilitate applications of mean field games,  we are motived to prove the required feedback regularity property  and provide verifiable conditions for it. This in turn motivates us to undertake the sensitivity analysis for HJB equations with respect to a functional parameter. Further,  we associate HJB equations to a class of general Markov processes, in order to provide a unified framework where the sensitivity analysis for general HJB equations with respect to functional parameters are studied. This is the main aim of this paper.

This work has two-fold contributions: 
first, this work presents a unified framework where the sensitivity analysis for general HJB equations with respect to a functional parameter is studied. For the stochastic control theory, these sensitivity results promote the understanding on how individual's optimal payoff depends on certain external observable functional parameters. These results have wide range of applications e.g. in financial markets,  economics and marine ecology. 
Second, for the mean field games theory, we prove the required feedback regularity property  and provide verifiable conditions for it. Our sensitivity results  enables applications of mean field games in solving real complex systems. 

The organisation of this paper is as follows. Section \ref{Problem description} gives the detailed description and the derivation of the problem in question. Section \ref{Preliminaries} recalls some basic but heavily used concepts for solving the problem in question. These concepts include operators, propagators and variational derivatives. The main results of this paper are presented in Section \ref{Main results}. Proofs for these results are collected in the appendix section. 
In Section \ref{Application}, as an application, we apply these sensitivity results to the study a mean field games model and give verifiable conditions for the feedback regularity property \eqref{FBR}.

\section{Problem description}
\label{Problem description}
In this section, we present the problem we aim to study, starting with the description of a stochastic control problem with a general Markov dynamics. 

Let $\C:=C_\infty(\R^d)$  be the Banach space of bounded continuous functions $f:\R^d\to \R$ with $\lim_{x\rightarrow\infty}f(x)=0$, equipped with norm $\|f\|_{\C}:=\sup_x|f(x)|$.
We shall denote by $\C^1:=C^1_\infty(\R^d)$ the Banach space of continuously differentiable and bounded functions $f:\R^d\to\R$  such that the derivative $f'$ belongs to $\C$, equipped with the norm $\|f\|_{\C^1}:=\sup_x|f(x)|+\sup_x |f'(x)|$, and by $\C^2:=C^2_\infty(\R^d)$ the Banach space of twice continuously differentiable and bounded functions $f:\R^d\to \R$ such that the first derivative $f'$ and the second derivative $f''$ belong to $\C$, equipped with the norm $\|f\|_{\C^2}:=\sup_x |f(x)|+\sup_x |f'(x)|+\sup_x |f''(x)|$. Let $\C_{Lip}:=C_{Lip}(\R^d)$ denote the space of Lipschitz continuous functions $f:\R^d\to \R$, equipped with the norm $\|f\|_{\C_{Lip}}:=\sup_x|f(x)|+\sup_{x,y}\frac{|f(x)-f(y)|}{|x-y|}$.

\begin{remark} 
Note that $\C^2=C^2_\infty(\R^d)$ is a Banach space which is densely and  continuously embedded in $\C=C_\infty(\R^d)$.
Depending on the modelling assumption,  the Banach space $\C^2$ can be replaced by other examples of functions spaces, such as the space of H{\"o}lder continuous functions, and our methods can be applied in a similar way.
\end{remark}
Let $T>0$ be fixed  and  $\U$ be a compact subset of a Euclidean space, interpreted as the set of admissible controls, with the Euclidean norm $|\cdot|$. Here (and in such like expressions), $\cdot$ denotes a dynamic variable. Take $\M$ to be a bounded, convex, closed subset of another Banach space $\S$, equipped with the norm  $\|\cdot\|_{\S}$. In applications, very often the Banach space $\S$ is taken as the dual space $(\C^2)^*$ of $\C^2$ and the set $\M$ is taken as the set of probability measures on $\R^d$, which is denoted by $\P(\R^d)$.

Now we introduce the stochastic dynamics associated to the stochastic control problem. Specifically, the dynamics is described by a family of bounded linear operators 
 \begin{equation}\label{A}
 \{A[t,\mu,u]:  t\in[0,T],\mu\in \M,u\in \U\}
 \end{equation}
where for all $(t,\mu,u)\in [0,T]\times \M\times \U$, $A[t,\mu,u]: \C^2\to \C$.
For each $(t,\mu,u)\in [0,T]\times \M\times \U$, the linear operator $A[t,\mu,u]$ is assumed to generate a Feller process with values in $\R^d$ and to be of the form
\begin{equation}
\label{fellergeneratorwithdriftcont}
A[t,\mu,u]f(z)=(h(t,z,\mu,u), \nabla f(z)) + L[t,\mu] f(z),
\end{equation}
where the coefficient $h: [0,T]\times\R^d\times\M\times \U\to \R^d$ is a vector-valued function.
For each pair $(t,\mu)\in [0,T]\times \M$, the linear bounded operator $L[t,\mu]:\C^2\to \C$ is of L\'evy-Khintchine form with variable coefficients:
{\small
\begin{equation}
\begin{split}
&L[t,\mu]f(z)=\frac{1}{2}(G(t,z,\mu)\nabla,\nabla)f(z)+ (b(t,z,\mu),\nabla f(z))\label{L}\\
 +&\int_{\R^d} (f(z+y)-f(z)-(\nabla f (z), y){\bf 1}_{B_1}(y))\nu (t,z,\mu,dy)
\end{split}
\end{equation}
}
where $\nabla$ denotes the gradient operator and $(G(t,z,\mu)\nabla,\nabla)=\sum _{i,j}G_{i,j} (t,z,\mu)\frac{\partial ^2}{\partial z_i\partial z_j}$; ${\bf 1}_{B_1}$ denotes the indicator function of the unit ball in $\R^d$ ;  for each $(t,z,\mu)\in [0,T]\times\R^d\times\M$, $G(t,z,\mu)$ is a symmetric non-negative matrix, $b(t,z,\mu)$ is a vector, $\nu(t,z,\mu,\cdot)$ is a L\'evy measure on $\R^d$, i.e.
\begin{equation}
\label{condlevy0}
\int_{\R^d} \min (1,|y|^2)\nu(t, z, \mu, dy) <\infty,
\end{equation}
with $\nu (t, z, \mu, \{0\})=0$. We assume that the mappings $(t,z,\mu)\to G(t,z,\mu)$, $(t,z,\mu)\to b(t,z,\mu)$ and $(t,z,\mu)\to \nu(t,z,\mu, \cdot)$ are Borel measurable with respect to the Borel $\sigma$-algebra in
$[0,T]\times\R^d\times\M$.

\begin{remark}
It worths noting on the special structure of operator $A$ in \eqref{fellergeneratorwithdriftcont}. The operator $A$ depends on three parameters: time $t$, control $u$ and unusually a $\mu\in\M$. Notice that the control parameter $u$ appears only in the drift coefficient $h$, but not in the operator $L$. In plain words, in the stochastic control problem, it is only allowed to control the deterministic component of the dynamics. Here the operator $L$ \eqref{L} models the noise, i.e. the stochastic component of the dynamics. The noise modeled by the L\'evy-Khintchine form operator $L$ in \eqref{L} includes Gaussian noise (see an example in Remark \ref{example of h and L} ) and L\'evy stable noise (see an example in \eqref{eqstable}).
\end{remark}

Denote by $\{u.\}=\{u_t\in \U: t\in[0,T]\}$ the control process and by $C([0,T], \M)$  the space of continuous  curves $\{\mu.\}=\{\mu_t\in \M$, $t\in [0,T]\}$. For given $\{u.\}$ and $\{\mu.\}\in C([0,T], \M)$, let $(X_t^{\{\mu.\},\{u.\}}: t\in[0,T])$ be a controlled stochastic process on a probability space $(\Omega, \FC, \PC)$ with values in $\R^d$ and generated by the family of operators $\{A[t,\mu,u]:  t\in[0,T],\mu\in \M,u\in \U\}$ in \eqref{A} of the form \eqref{fellergeneratorwithdriftcont}.
For notational brevity, in the following we write $(X_t: t\in[0,T])$ instead of $(X_t^{\{\mu.\},\{u.\}}: t\in[0,T])$.

\begin{remark} \label{example of h and L}
  \begin{enumerate}
  \item
   A simple example of the function $h$ in  \eqref{fellergeneratorwithdriftcont} is
    \[
    h(t,z,\mu,u)=\int_{\R^d} a(t,z, y,u)\mu(dy).
    \]
    with $a:  [0,T]\times\R^d\times\R^d\times \U\to \R^d$.
    \item
 In the special case of \eqref{L} where  $L[t,\mu]= \frac{1}{2}(G(t,z,\mu)\nabla,\nabla)$, the process associated with  the operators \eqref{A}  can also be described by  stochastic differential equations
\[
dX_t=h(t,X_t, \mu_t, u_t)\, dt+\sigma (t,X_t, \mu_t) dW_t
\]
with $G(t, z,\mu)=tr\{\sigma (t, z,\mu)\sigma^T (t, z,\mu)\}$, $W_t$ being a standard Brownian motion in $\R^d$.
\end{enumerate}
\end{remark}
Given the observable functional parameter $\{\mu.\}\in C([0,T], \M)$,  the objective is to maximize the expected total payoff
  \[
   \E \left[ \int_t^T J(s,X_s,\mu_s,u_s) \, ds +V^T (X_T,\mu_T)\right]
\]
at any time $t\in[0,T]$ over a suitable class of controls $\{u_t\in\U, t\in[0,T]\}$  with a running cost function $J: [0,T] \times \R^d\times \M\times \U \to \R$ and a terminal cost function $V^T: \R^d\times \M\to \R$.
Therefore, for a given parameter $\{\mu.\}\in C([0,T],\M)$, the value function $V:[0,T]\times \R^d\to \R$ starting at time $t$ and state $x$  is defined by
 \begin{align}\label{valuefunction}
&V(t,x; \{\mu.\}):=\\
&\sup_{\{u.\}} \E_x \left[ \int_t^T J(s,X_s,\mu_s,u_s) \, ds +V^T (X_T,\mu_T)\right].\notag
\end{align}

\begin{remark}
This type of stochastic control problems with a functional parameter not only arises from the study of mean field games (see more detailed discussion in \cite{KLY2012}  and Section \ref{Application} of this paper), but also appears in financial markets. For example, in the context of designing optimal trading strategies in high frequency trading based on limit order book dynamics, the functional parameter $\{\mu.\}$ can be interpreted as the bid/ask price distribution. 
\end{remark}

By standard arguments from dynamic programming principle and assuming appropriate regularity, the value function $V$ satisfies the Hamilton-Jacobi-Bellman (HJB) equation
\begin{align}\label{HJB_Mu}
\begin{split}
\displaystyle
-\frac{\partial{V}}{\partial{t}}(t,x;\{\mu.\}) &= H_{t}(x,\nabla V(t,x;\{\mu.\}), \mu_t) \\
&\hspace{1.5em}+L[t,\mu_t]V(t,x;\{\mu.\})\\
V(T,x;\{\mu.\})&=V^T(x;\mu_T),
\end{split}
\end{align}
 where the Hamiltonian $H:[0,T]\times \R^d\times\R^d\times \M \to \R$ is defined by
\begin{equation}\label{HJB H_Mu}
H_{t}(x,p, \mu)=\max_{u\in\U} ( h(t,x,\mu,u)p+J(t,x,\mu,u)).
\end{equation}

The aim of this paper is to investigate the sensitivity of the solution $V(\cdot,\cdot;\{\mu.\})$ to the HJB equation \eqref{HJB_Mu} with respect to the functional parameter $\{\mu.\}\in C([0,T], \M)$. 
In fact, we prove in Section \ref{Main results} that  the unique solution $V(\cdot,\cdot;\{\mu.\})$ to  the HJB equation \eqref{HJB_Mu} and its spatial gradient $\nabla V(\cdot,\cdot;\{\mu.\})$ are Lipschitz continuous  uniformly with respect to  $\{\mu.\}$.

\section{Preliminaries}\label{Preliminaries}

In this section, we recall some concepts which are used in the following of the paper.
Let $\B$ and $\D$ denote some Banach spaces. For a  function
$F:\D\to\B$, its {\it variational derivative} $D_\chi F(\mu)$ at $\mu\in \D$ in the direction $\chi\in \D$ is defined as
\[
D_\chi F (\mu)=\lim_{s\to 0} \frac{F(\mu+s\chi)-F(\mu)}{s}
\]
if the limit exists. F is said to be differentiable at $\mu\in\D$ if the limit exits for all $\chi\in\D$.  At each point $\mu\in \D$, the derivative  defines a function $D_. F(\mu):\D\to\B$.

Let ${\mathcal L}(\D,\B)$ denote the space of linear bounded operators from $\D$ to $\B$ and it is equipped with the usual operator norm
$\|\cdot\|_{\D\to\B}$.

For the analysis of time non-homogeneous evolutions, we need the notion of a propagator. A family of mappings $\{U^{t,r}\}$ from $\B$ to $\B$, parametrized by the pairs of numbers $r\leq t$ (resp. $t\leq r$) is called a {\it (forward) propagator} (resp. a {\it backward propagator}) in $\B$, if $U^{t,t}$ is the identity operator in $\B$ for all $t\geq 0$ and the following {\it chain rule}, or {\it propagator equation}, holds for $r\leq s\leq t$ (resp. for $t\leq s\leq r$):
$$U^{t,s}U^{s,r}=U^{t,r}.$$
Sometimes, the family $\{U^{t,r}, t\leq r\}$ is also called a {\it two-parameter semigroup}.

A backward propagator $\{U^{t,r}, t\leq r\}$ of bounded linear operators on the Banach space $\B$ is called {\it strongly continuous} if for every $f\in\B$, the mappings
$$
t\mapsto U^{t,r}f \,\text{for all }t\leq r,\,\text{ and }\,
r\mapsto U^{t,r}f \,\text{for all }t\leq r,
$$
are continuous as mappings from $\R$ to $\B$. By the principle of uniform boundedness if $\{U^{t,r}, t\leq r\}$ is a strongly continuous propagator of bounded linear operators, then the norms of $\{U^{t,r}, t\leq r\}$ are uniformly bounded  for $t,r$ in any compact interval.

Assume that the Banach space $\D$ is a dense subset of $\B$ and continuously embedded in $\B$. Suppose $\{U^{t,r}, t\leq r\}$ is a strongly continuous backward propagator of bounded linear operators on a Banach space $\B$ with the common invariant domain $\D\subset \B$, i.e. if $f\in\D$ then $U^{t,r}f\in \D$ for all $t\leq r$. Let $\{L_t, t\geq0\}$ be a family of operators $L_t\in {\mathcal L}(\D,\B)$, depending continuously on $t$. The family $\{L_t, t\geq0\}$ is said to {\it generate} $\{U^{t,r}, t\leq r\}$ on $\D$ if,
for any $f\in \D$, we have, for all $t\leq s\leq r$
\begin{align}
\label{Generates}
 \frac{d}{ds}U^{t,s}f = U^{t,s}L_sf, \quad \frac{d}{ds}U^{s,r}f = -L_sU^{s,r}f.
\end{align}
The derivatives exist in the norm topology of $\B$ and if  $s=t$ (resp. $s=r$) they are assumed to be only a right (resp. left) derivative.

One often needs to estimate the difference of two propagators when the difference of their generators is available. To this end, we shall often use the following rather standard trick.
\begin{prop}
\label{prop-propergatorProperty}
For $i=1,2$ let $\{L_t^i, t\geq0\}$ be a family of operators $L_t^i\in {\mathcal L}(\D, \B)$, depending continuously on $t$, which generates a backward propagator
$\{U_i^{t,r}, t\leq r\}$ in $\B$ satisfying
$$
a_1:=\sup_{t\leq r}\max\left\{\|U_1^{t,r}\|_{\B\to \B}, \|U_2^{t,r}\|_{\B\to \B} \right\}<\infty.
$$
If $\D$ is invariant under $\{U_1^{t,r}, t\leq r\}$ and
$$
a_2:=\sup_{t\leq r}\|U_1^{t,r}\|_{\D\to \D}<\infty,
$$
then, for each $t\leq r$,
\begin{equation}
\label{PropergatorProperty }
U_2^{t,r}-U_1^{t,r}=\int_t^rU_2^{t,s}(L^2_s-L_s^1)U_1^{s,r}ds
\end{equation}
and
\begin{align}
&\|U_2^{t,r}-U_1^{t,r}\|_{\D\to \B} \label{PropergatorPropertya}\\
\le &a_1a_2(r-t) \sup_{t\leq s\le r}\|L^2_s-L_s^1 \|_{\D\to \B}.\notag
\end{align}
\end{prop}
\proof
By (\ref{Generates}), we get
\begin{equation*}
\begin{split}
U_2^{t,r}-U_1^{t,r}&=U_2^{t,s}U_1^{s,r}\big|_{s=t}^r=\int_t^r\frac{d}{ds} \left(U_2^{t,s}U_1^{s,r} \right) ds\\
&=\int_t^r U_2^{t,s} L_s^2 U_1^{s,r}-U_2^{t,s} L_s^1 U_1^{s,r}ds\\
&=\int_t^rU_2^{t,s}(L^2_s-L_s^1)U_1^{s,r}ds,
\end{split}
\end{equation*}
which implies both \eqref{PropergatorProperty } and \eqref{PropergatorPropertya}.
\qed

\section{Main results}\label{Main results}

In this section, first we  list the assumptions needed for the discussion. Then based on these assumptions, we show that for each fixed parameter $\{\mu.\}\in C([0,T], \M)$, the HJB equation \eqref{HJB_Mu} is well posed and we prove the Lipchitz continues dependence of the solution to \eqref{HJB_Mu} with respect to the functional parameter $\{\mu.\}\in C([0,T], \M)$. 

\subsection{Main assumptions}\label{assumptions}

For any $\mu\in\M$, define the set $\M - \mu:= \{\eta-\mu: \eta\in\M\}$, which, as a subset of $\S$, is equipped with the norm $\|\cdot\|_\S$. In the analysis below,we need the following assumptions:

{\bf (A1)}: the Hamiltonian $H:[0,T]\times \R^d\times\R^d\times \M \to \R$ is continuous in $t$  and   Lipschitz continuous uniformly in $x$ on bounded subsets of $p$. Furthermore, it is
 Lipschitz continuous uniformly in $p$, that is there exists a constant $c_1$ such that
 for all $x\in \mathbf{R}^d$, $\mu \in  \M$ and $t\in[0,T]$ we have for $p,p'\in \mathbf{R}^d$
\begin{equation}
\label{eq1thweak_existencem}
|H_{t}(x,p,\mu)-H_{t}(x,p',\mu)|\leq c_1|p-p'|.
\end{equation}
It is bounded in $p=0$, that is there  exists a constant $c_2>0$ such that for all $x\in \mathbf{R}^d,\, \mu \in  \M,\, t\in[0,T]$
\begin{equation}
\label{eq2thweak_existencem}
|H_{t}(x,0, \mu)|\leq c_2.
\end{equation}
For each $t\in [0,T]$ and $x,p\in\R^d$, the function $\mu\mapsto H_{t}(x,p, \mu)$ is differentiable in any direction $\chi\in\M - \mu$, such that $(t,x,p,\mu)\mapsto D_\chi H_{t}(x,p, \mu)$  is continuous  and satisfies that for each bounded set $B\subset \R^d$ there exists a constant $c_3>0$ such that for all $t\in [0,T], x\in \R^d, \mu\in \M$
\begin{equation}
\label{eqassumonder1}
\sup_{p\in B}|D_\chi H_{t}(x,p, \mu)|\le c_3\|\chi\|_{\S} .
\end{equation}

{\bf (A2)}:  the mapping
\begin{align*}
 [0,T]\times \M\to {\mathcal L}(\C^2,\C), \qquad (t,\mu)\mapsto L[t,\mu]
\end{align*}
is continuous. For any $\{\mu.\} \in C([0,T], \M)$, the operator curve $\{L[t,\mu_t]: t\in[0,T]\}$ generates a strongly continuous backward propagator $\{U^{t,s}_{\{\mu.\}},t\leq s\}$ of operators $U^{t,s}_{\{\mu.\}}\in {\mathcal L}(\C,\C)$ with the common invariant domains $\C^2$
and $\C^1$. There exists a constant $c_4>0$ such that
 for all $0\leq t\leq s\leq T$ we have
\begin{align}\label{BDD 2m}
\max &\left\{\|U^{t,s}_{\{\mu.\}}\|_{\C\to \C}, \|U^{t,s}_{\{\mu.\}}\|_{{\C^1}\to{\C^1}}, \|U^{t,s}_{\{\mu.\}}\|_{\C^2 \to\C^2}\right\}\notag\\
&\leq c_4.
	\end{align}
The propagator has a {\it smoothing property}, that is
for each $0\leq t<s\leq T$ we have
 \begin{equation}\label{smoothingm}
U^{t,s}_{\{\mu.\}}: \C\to \C^1, \quad
U^{t,s}_{\{\mu.\}}: \C_{Lip}\to \C^2,
\end{equation}
and there exists a $\be \in (0,1)$ and constants $c_5,c_6>0$ such that
\begin{align}
&\|{U}^{t,s}_{\{\mu.\}}\phi\|_{\C^1} \leq c_5(s-t)^{-\be}\|\phi\|_{\C},\label{smooth property 2m-1}\\
&\|{U}^{t,s}_{\{\mu.\}}\psi\|_{\C^2} \leq c_6(s-t)^{-\be} \|\psi\|_{\C_{Lip}}\label{smooth property 2m}
\end{align}
for all $\phi\in \C$ and $\psi\in \C_{Lip}$.

(ii) for any $t\in[0,T]$, the mapping $\mu\mapsto L[t,\mu]$ is  differentiable in any direction $\chi\in\M-\mu$, such that the mapping
$\mu\mapsto D_\chi L[t,\mu]$
is continuous. There exists a constant $c_7>0$ such that for each $\mu\in \M$ and $\chi\in \M-\mu$ we have for all $t\in[0,T]$
\begin{equation}
\label{eqassumonder2}
\left\|D_\chi L[t,\mu]\right\|_{\C^2 \to\C}\leq c_7 \|\chi \|_{\S};
\end{equation}

{\bf (A3)}: for any $\mu \in  \M$, the mapping $x\mapsto V^T(x; \mu)$ is twice continuously differentiable, and for each $x\in\R^d$ the  mapping $\mu\mapsto V^T(x; \mu)$ is  differentiable in any direction $\chi\in\M-\mu$ such that
the mapping $(x,\mu)\mapsto  D_\chi V^T(x;\mu)$ is continuous. There exists a constant
$c_{8}>0$ such that
\begin{equation}
\label{eqassumonder3}
\|D_\chi V^T(\cdot;\mu)\|_{\C^1}\le c_8 \|\chi\|_{\S}.
\end{equation}

\begin{remark}
\begin{enumerate}
\item
If the Banach space $\S$ is given as the Euclidean space $\R$, then  $D_{\chi}$ corresponds to the standard partial derivatives and are denoted by $\partial/\partial \al$
for $\al\in\R$.
\item
The assumptions {\bf (A1)} {\bf (A2)}  and {\bf (A3)}  are stated in an abstract setting. These assumptions are made concrete in Example \ref{example assumptions} in Section \ref{Application} in an application to mean field games. 
\end{enumerate}
\end{remark}

The smoothing conditions \eqref{smoothingm}, \eqref{smooth property 2m-1} and \eqref{smooth property 2m} in assumption {\bf (A2)} are essential and critical in the following analysis. Let us show two basic examples which satisfy assumption {\bf (A2)}:
the diffusion operator
{\small
\begin{align}
&L[t,\mu]f(x)\notag\\
=&\frac{1}{2}(\si ^2(t,x,\mu) \nabla,\nabla)f(x)+(b(t,x,\mu), \nabla f)(x)\label{eqdifpr}
\end{align}
}
with smooth enough functions $b,\si$, see e.g. in \cite{PorEid84} and references therein.  The operators
$\{L[t,\mu],t\in[0,T]\}$ generate the stochastic process $(X(t),t\in[0,T])$ which obeys
 the stochastic differential equation
\[
dX_t = b(t, X_t, \mu_t)dt + \sigma(t, X_t, \mu_t)dW_t,
\]
where $W$ is a standard Brownian motion. 

Another example is given by  stable-like processes with the generating family
\begin{align}
&L[t,\mu]f(x)\notag\\
=&a(t,x) |\De|^{\al (x)/2}+(b(t,x,\mu), \nabla f)(x)\label{eqstable}
\end{align}
with smooth enough functions $a, \al$ such that the range of $a$ is a compact interval of positive numbers and the range of $\al$ is a compact subinterval of $(1,2)$.

In both cases, each operator $U_\mu^{t,s}$, $t\leq s$,  has a kernel, e.g. it is given by
\begin{equation}
U_\mu^{t,s} f(x)= \int G_\mu(t,s,x,y) f(y) dy
\end{equation}
with a certain Green's function $G_\mu$, such that
for every $x\in\R^d$ and $t\leq s$,
\begin{equation}
\sup_{\mu\in\M}\int_{\R^d} |\nabla_x G_\mu(t,s,x,y) | dy\leq  c(s-t)^{-\be}
\end{equation}
for a constant  $c>0$. For the case of standard Brownian motion (i.e. $L[t,\mu]f(x)=\frac{1}{2}\Delta f(x)$), we have explicit expression with $\beta=\frac{1}{2}$; for general diffusion processes \eqref{eqdifpr}, we know the function $G$ is bounded by Gaussian kernel  and its derivate of $G$ is bounded by the one of the Gaussian kernel  with $\beta=\tfrac{1}{2}$; for the stable process  \eqref{eqstable}, we know  the function $G$ is bounded by stable density and its derivate of $G$ is bounded by the one of the stable density with $\beta=(\inf_x \al(x))^{-1}$. Then \eqref{smooth property 2m-1}  is implied. \eqref{smooth property 2m} can be obtained in the similar but lengthy manner by differentiation with respect to $x$, see details in \cite{VK2} and references therein.

\subsection{Well-posedness of HJB equations}
In this subsection, we prove the well-posedness of the HJB equation \eqref{HJB_Mu} for any given $\{\mu_t:t\in [0,T]\}\in C([0,T],\M)$.
For this purpose,  we generalise the explicit dependence of the functions $H, L, V^T$  on the parameter $\mu\in\M$, and  encode their dependence on $\mu_t$ through time $t$. Specifically, we consider the  general Cauchy problem
\begin{align}\label{G_HJB}
\begin{split}
-\frac{\partial{V}}{\partial{t}}(t,x) &= H_t(x,\nabla V(t,x)) +L_tV(t,x)\\
V(T,x)&=V^T(x)
\end{split}
\end{align}
with the Hamiltonian $H:[0,T]\times\R^d\times\R^d\to \R$ defined by
\begin{equation}\label{HJB 2b}
H_t(x,p)=\max_{u\in\U} ( h(t,x,u)p+J(t,x,u)).
\end{equation}
and the operator $L_t: \C^2\to \C$ for each $t\in [0,T]$. Notice that the wellposedness of \eqref{G_HJB} immediately implies the wellposedness of \eqref{HJB_Mu} for each fixed 
$\{\mu_t:t\in [0,T]\}\in C([0,T],\M)$
under the same conditions.

By Duhamel's principle (c.f. \cite{E2010}), if $V$ is a classical solution of \eqref{G_HJB}, then $V$ is also
a {\it mild solution} of \eqref{G_HJB}, i.e. it satisfies
{\small
\begin{align}
V(t,x) = &(U^{t,T} V^T(\cdot))(x) \notag\\
&+ \int_t^T U^{t,s}H_s(\,\cdot\,,\nabla V(s,\cdot))(x) ds\label{mild_value_function}
\end{align}
}
for all $t\in [0,T]$ and $x\in\R^d$.

For the sensitivity analysis of this work, it is sufficient to consider only a mild solution. For this reason, we will establish
the existence of a unique mild solution. In this subsection, we mostly follow Chapter 7 in \cite{VK2}, where one also can find details  for existence of a classical solution. We present this result for completeness on the level of generality which is required by what follows.

\begin{theorem}
\label{weak_existence}
Assume conditions {\bf (A1)} and {\bf (A2)}. If the terminal data $V^T(\cdot)$ is in $\C^1$, then there exists a unique mild solution $V$ of \eqref{G_HJB}, satisfying $V(t,\cdot)\in \C^1$ for all $t\in[0,T]$.
\end{theorem}
\proof see the proof in the appendix \ref{Proof to the Theorem weak_existence }.
\qed

By the wellposedness of equation \eqref{mild_value_function}, its solution defines a propagator
in $\C^1$.  The regularising term $L_tV(t,x)$ in the HJB equation \eqref{G_HJB} allows us to get a more regular solution  than a mild solution.  Standard arguments, see e.g. \cite{FlSo}, show
that this mild solution is a viscosity solution to the original equation \eqref{G_HJB}
and it solves the corresponding optimization problem.

\subsection{Sensitivity analysis for HJB equations}\label{Sensitivity analysis of HJB }

In this subsection, we  analyse the dependency of the solution of the HJB equation \eqref{HJB_Mu} on the functional parameter $\{\mu.\}\in C([0,T],\M)$. Under the conditions {\bf (A1)} and {\bf (A2)}, Theorem \ref{weak_existence} guarantees the existence of a unique mild solution $V(\cdot,\cdot; \{\mu.\})$ of \eqref{HJB_Mu} for each fixed  $\{\mu.\}\in C([0,T],\M)$.

The following observation plays an important role in this work.
Let  $\{\mu^1_.\},\{\mu^2_.\}$ be in $C([0,T], \M)$ and let $\alpha \in [0,1]$.
Since $\M$ is convex, the curve
$$\{\mu^1_.\}+\alpha \{(\mu^2-\mu^1)_.\}:=\{\mu^1_t+\alpha (\mu^2_t-\mu^1_t), t\in[0,T]\}$$
belongs to $C([0,T],\M)$. Thus, we can define the function $V_\alpha:[0,T]\times\R^d\to \R$
\begin{equation}\label{Valpha}
V_\alpha(t,x):=V\big(t,x; \{\mu^1_.\}+\alpha \{(\mu^2-\mu^1)_.\}\big).
\end{equation}
and have the relation
\begin{align}
&V(t,x; \{\mu^2_.\})-V(t,x;  \{\mu^1_.\})\notag\\
=&V_{1}(t,x)-V_{0}(t,x).\label{estimate}
\end{align}
Furthermore, for $\{\mu^1_.\},\{\mu^2_.\}\in C([0,T], \M)$ define
\begin{align}\label{Halpha}
H_{\al,t}(x,p)
:&=H_{t}(x,p, \mu^1_t+\alpha (\mu^2_t-\mu^1_t))
 \\
&=\max_{u\in\U} ( h(t,x,\mu^1_t+\alpha (\mu^2_t-\mu^1_t),u)p\notag\\
&\hspace{1em}+J(t,x,\mu^1_t+\alpha (\mu^2_t-\mu^1_t),u))\notag
\end{align}
\begin{equation}\label{Lalpha}
\begin{split}
L_{\al}[t]:=L(t, \mu^1_t+\alpha (\mu^2_t-\mu^1_t))\\
\end{split}
\end{equation}\begin{equation}\label{VTalpha}
\begin{split}
V_{\al}^T(x):=V^T(x;\mu^1_T+\alpha (\mu^2_T-\mu^1_T))
\end{split}
\end{equation}
 with $\al\in[0,1]$, $t\in[0,T]$ and $(x,p)\in\R^{d}\times\R^d$. Then the sensitivity analysis of the solution of \eqref{HJB_Mu} with respect to a function parameter $\{\mu.\}\in C([0,T],\M)$ can be reduced to the one of the solution to the following Cauchy problem with respect to a real parameter $\alpha\in[0,1]$:
\begin{align}
\frac{\partial{V_\alpha}}{\partial{t}}(t,x)& =-H_{\al,t}(x, \nabla V_\alpha(t,x))- L_{\al}[t]V_\alpha (t,x)\notag\\
V_\alpha(T,x)&=V_\alpha^T(x).\label{V_alpha}
\end{align}

The sensitivity analysis with respect to $\al\in [0,1]$ consists of two steps. In the first step, we omit the Hamiltonian term in \eqref{V_alpha} and only consider the sensitivity of the evolution $V_{\al}(t,\cdot)=U^{t,T}_{\al} V^T_{\al}(\cdot)$,  where for each $\alpha\in [0,1]$, the propagator  $\{U_{\al}^{t,s}:t\leq s\}$ is generated by the family of operators $\{L_{\al}[t]: t\in[0,T]\}$. This step  not only serves as an intermediate step towards our full scheme, but also unveils some interesting observations, including the formula \eqref{smooth property 6} for the differentiation of an operator with respect to a parameter. In the second step,  we include the Hamiltonian term and complete the analysis.

\begin{theorem} \label{Sensitivity 1}
Assume conditions {\bf (A2)} and {\bf (A3)}. Define $W:[0,1]\times [0,T]\times \R^d\to\R^d$
\begin{equation}\label{smooth property 3}
W_{\al}(t,x)=U^{t,T}_{\al} V^T_{\al}(x).
\end{equation}
Then for each $t\in [0,T]$ and $x\in\R^d$, the mapping $\alpha\to W_\al(t,x)$ is Lipschitz continuous  with uniformly bounded Lipschitz constants, more precisely for every $\al_1, \al_2\in[0,1]$ with $\al_1\neq \al_2$, there exists a constant $c>0$ such that for every $t\in[0,T]$
{\small
\begin{align*}
&\frac{\|W_{\al_1}(t,\cdot)-W_{\al_2}(t,\cdot)\|_{\C^1}}{|\al_1-\al_2|}\notag\\
\leq&
c(T-t)^{1-\beta}
\underset{\begin{subarray}{c}
  s\in[t,T] \\
  \gamma\in[\al_1,\al_2]
  \end{subarray}}{\sup}\left\| \frac{\partial L_\gamma}{\partial \al}[s]\right\|_{\C^2\to\C}\|V^T_{\al_2}(\cdot)\|_{\C^2}\\
  &+ c \sup_{\gamma\in[\al_1,\al_2]}\left\|\frac{\partial V_\gamma^T}{\partial \al}\right\|_{\C^1}.
\end{align*}
}
\end{theorem}
\proof
See the proof in the appendix \ref{Proof to the Theorem Sensitivity 1}.
\qed

\begin{remark}
To clarify the notations $\alpha$ and $\gamma$ used in the derivative $\partial {V_\gamma^T}/{\partial \alpha}$: $V^T_\alpha$ denotes that the $\alpha$ is the variable of the function $V^T$ and $\frac{\partial}{\partial \alpha} $ denotes the derivative of a function with respect to the variable $\alpha$. The $\gamma$ is the value of the variable  used to calculate the value of the derivative. $\frac{\partial V_\gamma^T}{\partial \alpha}$ denotes the value of the derivative of the function $V^T$ with respect to $\alpha$ at the point $\alpha =\gamma$.
\end{remark}

In this work, we are only concerned with the  Lipschitz continuity of the solution of the HJB with respect to the parameter. However, it is also interesting to know whether the mapping $\al\mapsto W_{\al}(t,\cdot)$ is differentiable for each $t\in [0,T]$. For the completeness, the next proposition will show  the existence of the derivative  $\frac{\partial W_{\al}}{\partial \al}(t,\cdot)$ in $\C$ and present its' explicit expression.
\begin{prop}
\label{prop}
Assume conditions {\bf (A2)} and {\bf (A3)}. Then
\begin{enumerate}
\item[{\rm (i)}] for each $0\leq t<s\leq T$, the mapping $\al\mapsto U_\alpha^{t,s}$ is differentiable and  the derivative $\frac{\partial U_\alpha^{t,s}}{\partial \alpha }$  has the representation
\begin{equation}
\label{smooth property 6}
\frac{\partial U_\alpha^{t,s}}{\partial \alpha}
=\int_t^s {U}^{t,r}_{\alpha} \frac{\partial{L_\alpha}}{\partial{\alpha}}[r] U_\al^{r,s} \, dr.
\end{equation}

\item[{\rm (ii)}] for each $t\in[0,T]$, the mapping $\al\mapsto W_\al(t,\cdot)$ defined in \eqref{smooth property 3}  is differentiable as a function from $[0,1]$ to $\C$ and the partial derivative $\tfrac{\partial W_\alpha}{\partial \alpha}(t,\cdot)$
can be represented by
\begin{equation}
\label{smooth property 5}
\frac{\partial W_\alpha}{\partial \alpha} (t,\cdot) = {U}_\alpha^{t,T} \frac{\partial V^T_\alpha}{\partial \alpha} (\cdot)
+\frac{\partial U_\alpha^{t,T}}{\partial \alpha}V^T_\alpha(\cdot).
\end{equation}
\end{enumerate}
\end{prop}
\proof
See the proof in the appendix \ref{Proof for the Proposition prop}.
\qed

\begin{remark}
If one would have that, for each $0\leq t<s\leq T$, the mapping $\al\mapsto U_{\al}^{t,s}$ is continuous from $[0,1]$ to $\mathcal{L}(\C, \C^1)$, then $\lim_{\al_1\to \al}U_{\al}^{t,s}$ exists as an operator from $\C$ to  $\C^1$. Then one would have  that the mapping $\al\mapsto W_\al(t,\cdot)$ is differentiable as a function from $[0,1]$ to $\C^1$.
\end{remark}

\begin{theorem} \label{Sensitivity 2}
Assume the conditions {\bf (A1), (A2)} and {\bf (A3)}. Then the following statements hold:
\begin{enumerate}
\item[{\rm (a)}]
For any $T>0$, the mild solution $V_\alpha$ to \eqref{V_alpha} is Lipschitz continuous with
respect to $\alpha$ i.e. there exists a constant $c=c(T)>0$ such that for each $\al_1,\al_2\in[0,1]$ with $\al_1\neq \al_2$,
{\small
\begin{align}
&\sup_{t\in[0,T]}\frac{\| V_{\al_1}(t,\cdot)-V_{\al_2}(t,\cdot)\|_{\C^1}}{ |\al_1-\al_2|}\label{T10}\\[.4em]
&\hspace{-1cm} \leq c\Bigg( \sup_{\gamma\in[\al_1,\al_2]}\left\|\frac{\partial{V^T_\gamma}}{\partial{\alpha}}(\cdot)\right\|_{\C^1}
+\underset{\begin{subarray}{c}
  (t,p)\in\O \notag\\
  \gamma\in[\al_1,\al_2]
  \end{subarray}}{\sup}\left\|\frac{\partial H_{\gamma,t}(\cdot,p)}{\partial \alpha}\right\|_{\C}\notag\\[.4em]
&\hspace{-1cm}\hspace{1em}+ \underset{\begin{subarray}{c}
  t\in[0,T] \notag\\
  \gamma\in[\al_1,\al_2]
  \end{subarray}}{\sup}\left\|\frac{\partial L_\gamma}{\partial \alpha}[t]\right\|_{\C^2\to \C}  \left(\left\|V_{\al_2}^T(\cdot)\right\|_{\C^2}+1\right)\Bigg),\notag
\end{align}
}
where $\O=\{(t,p): t\in[0,T],|p|\leq \sup_{t\in[0,T]}\|V_\al(t,\cdot)\|_{\C^1}\}$.

\item[{\rm (b)}]   The mild solution $V$ to \eqref{G_HJB} and its spacial derivative $\nabla V$ are Lipschitz continuous uniformly with respect to $\{\mu.\}$, that is, for each $\{\mu^1_.\},\{\mu^2_.\}\in C([0,T], \M)$, there exists a constant $k>0$ such that
\begin{align}\label{V Lip}
&\sup_{t\in[0,T]}\|V(t,\cdot; \{\mu^1_.\})-V(t,\cdot; \{\mu^2_. \})\|_{\C^1}\notag\\
&\leq k \sup_{t\in [0,T]}\|\mu_t^1-\mu_t^2\|_{\S}
\end{align}
and
\begin{align}\label{V Lip_2}
&\sup_{t\in[0,T]}\|\nabla V(t,\cdot;  \{\mu^1_.\})-\nabla V(t,\cdot;  \{\mu^2_.\})\|_{\C}\notag\\
&\leq k \sup_{t\in [0,T]}\|\mu_t^1-\mu_t^2\|_{\S}.
\end{align}
\end{enumerate}
\end{theorem}
\proof
See the proof in the appendix \ref{Proof for the Theorem Sensitivity 2}.
\qed

\section{An application to mean field games}\label{Application}

In this section, we apply the sensitivity results in Theorem \ref{Sensitivity 2} to a mean field games model and to provide verifiable conditions for the feedback regularity condition.

Consider a continuous time dynamic game with a continuum of players and a terminal time $T>0$. By saying a continuum of players, we mean that all players are identical so the game is symmetric with respect to permutation of the players. Choose one of the players and call it the {\it reference player}. 

Take $\S=(\C^2)^*$ as the dual Banach space of $\C^2$ and $\M= \P(\R^d)$. Let $\mu_t\in \P(\R^d)$ denote the probability  distribution of the continuum of players in the state space $\R^d$ at the time $t$ and $\{\mu_t\in \P(\R^d):t\in[0,T]\}$ denote the (probability distribution) {\it measure flow}.
The controlled state dynamics of the reference player is modelled by a controlled stochastic process $(X_t:t\in[0,T])$ associated to the family of operators $A$ in \eqref{A}. At each time $t\in[0,T]$, the reference player knows only his own position $X_t$ and the distribution of the continuum of players $\mu_t\in \P(\R^d)$.
\begin{remark}
For a better understanding of the stochastic process $(X_t: t\in [0,T])$, one may think that  the controlled dynamics is described by a stochastic differential equation. For a very particular case of our model, set $L[t,\mu]= \frac{1}{2}(\sigma\nabla,\nabla)$ with a constant $\sigma$, i.e. the operator $L[t,\mu]$ generates a Brownian motion $\{\sigma W_t: t\geq 0\}$. Then one can write a stochastic differential equation corresponding to the generator \eqref{fellergeneratorwithdriftcont} as
\[
dX_t=h(t,X_t, \mu_t, u_t)\, dt+\sigma dW_t \quad \text{for all}\,\, t\geq 0.
\]
In fact, this is exactly the case which was considered in the initial work on the mean field games  \cite{HMC05, HCM3, HCM07, LL2007}.  In our framework, this controlled dynamics of each player is extended to an arbitrary Markov process with a generator \eqref{fellergeneratorwithdriftcont} depending on a probability measures $\mu$.
\end{remark}

The measure flow of the continuum of players in the state space $\R^d$, $\{\mu_t\in \P(\R^d):t\in[0,T]\}$, is the solution to the evolution equation
\begin{align}\label{limiting KE}
& \frac{d}{dt} \int_{\R^d} g(y)\,\,\mu_t(dy)\notag\\
& =\int _{\R^d} \left(A[t,\mu_t, u_t]g(y)\right)\,\, \mu_t(dy)
\end{align}
for a test function $g\in \C^2$, with a given initial value $\mu_0\in\P(\R^d)$.
The equation \eqref{limiting KE} is a controlled version of a {\it general kinetic equation} in weak form  and very often written in the compact form
\begin{equation}\label{dynamic-abs}
{\frac{d} {dt}} (g, \mu_t) =(A[t,\mu_t, u_t ]g, \mu_t).
\end{equation}
See \cite{KLY2012, KY2012} for more discussion on equation \eqref{dynamic-abs} and its well posedness with open-loop controls under rather general technical assumptions.  Let $C_{\mu_0}([0,T], \P(\R^d))$ be the set of continuous functions $t \rightarrow \mu_t$  with $\mu_t\in \P(\R^d)$ for each $t\in[0,T]$ and with the norm
\begin{align}
\label{D*}
  \|\mu\|_{(\C^2)^*}:=&\sup_{\|g\|_{\C^2}\leq 1}|(g,\mu)|\notag\\
  =&\sup_{\|g\|_{\C^2}\leq 1}\left|\int_{\R^d}g(x)\mu(dx)\right |.
\end{align}
In this game, the reference player faces an optimal control problem described by the HJB equation \eqref{HJB_Mu}. If the max is achieved only at one point, i.e. for any $(t,x,\mu, p)\in[0,T]\times\R^d\times\P(\R^d)\times\R^d$,
$$\arg\max_{u\in\U} ( h(t,x,\mu,u)p+J(t,x,\mu,u))$$
is a singleton, then one can derive the unique optimal control strategy from the solution to \eqref{HJB_Mu}. For any given measure flow $\{\mu_t: t\in [0,T]\}\in C_{\mu_0}([0,T], \P(\R^d))$, let the resulting unique optimal control strategy be denoted by
\begin{equation}\label{feedbacklaw}
\hat u(t,x;\{\mu_s: s\in [t,T]\})
\end{equation}
for all $t\in[0,T], x\in \R^d$. Substituting the optimal feedback control strategy \eqref{feedbacklaw} into \eqref{dynamic-abs} yields the closed-loop evolution equation for the distributions $\mu_t$
\begin{align}
{\frac{d}  {dt}} (g, \mu_t)& =(A[t,\mu_t, \hat u(t,x;\{\mu_s: s\in [t,T]\})]g, \mu_t)\label{Kinetic equation_coupled}
 \end{align}
The mean field game methodology amounts to find a value function $V$ for the representative player and a measure flow $\{\mu_.\}$ such that the two coupled equations \eqref{Kinetic equation_coupled}
with initial data $ \mu|_{t=0}=\mu_0$ and
\begin{equation}
-\frac{\partial{V}}{\partial{t}}(t,x) = H_{t}(x,\nabla V(t,x), \mu_t)+L[t,\mu_t]V(t,x)\label{HJB_coupled}
\end{equation}
with terminal data $V(T,\cdot;\mu_T)=V^T(\cdot;\mu_T)$ and 
\begin{equation}
H_{t}(x,p, \mu)=\max_{u\in\U} ( h(t,x,\mu,u)p+J(t,x,\mu,u))\label{singleton control}
\end{equation}
are satisfied at the same time, with the optimal control function $\hat u$ being the argmax in \eqref{HJB_coupled}. The HJB equation \eqref{HJB_coupled} is exactly the HJB equation \eqref{HJB_Mu}, which is the main object of this paper.

Since the controlled kinetic equation \eqref{Kinetic equation_coupled} is forward and the HJB equation \eqref{HJB_coupled} is backward, this system of coupled equations is referred to as a {\it coupled backward-forward system}. See \cite{KLY2012} for the full discussion on solving the coupled backward-forward system of equations \eqref{Kinetic equation_coupled} and \eqref{HJB_coupled}.

To solve this coupled backward-forward system \eqref{Kinetic equation_coupled}-\eqref{HJB_coupled}, it is critical that the resulting control mapping $\hat u$ \eqref{feedbacklaw} satisfies the so-called {\it feedback regularity} condition (see e.g. \cite{HCM3}), i.e. for any $\{\eta_t: t\in[0,T]\}$, $\{\xi_t: t\in[0,T]\}\in C_{\mu_0}([0,T], \P(\R^d))$, there exists a constant $k_1>0$ such that
\begin{align}
&\sup_{(t,x)\in [0,T]\times \R^d} \Big | \hat u(t,x ;\{\eta_s: s\in [t,T]\})\label{FBR}\notag\\
&\hspace{7em}-\hat u(t,x;\{\xi_s: s\in [t,T]\})\Big|\notag\\
&\leq   k_1\sup_{s\in [0,T]}||\eta_s -\xi_s||_{(\C^2)^*}.
\end{align}

\begin{theorem}
\label{thfeedbackHJB}

Suppose  $L$, $H$ and $V^T$ in \eqref{HJB_coupled}  satisfy the conditions {\bf (A1)}, {\bf (A2)}, {\bf (A3)} respectively. 
Assume additionally that 
the max in \eqref{singleton control} is achieved only at one point, i.e. for any $(t,x,\mu, p)\in[0,T]\times\R^d\times\P(\R^d)\times\R^d$
\begin{equation}\label{uni_u}
\arg\max_{u\in\U} ( h(t,x,\mu,u)p+J(t,x,\mu,u))
\end{equation}
is a singleton and the resulting control as a function of $(t,x,\mu, p)$ is continuous in $t\in[0,T]$ and Lipschitz continuous in $(x,\mu, p)\in\R^d\times\P(\R^d)\times\R^d$ uniformly with respect to $t$, $x$, $\mu$ and bounded $p$.
Then, the optimal feedback control strategy
$$\hat u (t,x;\{\mu_s: s\in [t,T]\})$$
derived via the HJB equation \eqref{HJB_coupled}, has the feedback regularity property defined in \eqref{FBR}.
\end{theorem}

\proof
By Theorem \ref{Sensitivity 2}, together with the assumption that the resulting unique control mapping is Lipschitz continuous in $(x,\mu,p)$, we conclude that the unique point of maximum in
the expression
\[
\max_{u\in\U} \{  h(t,x,\mu_t,u)\nabla V(t,x;\{\mu.\})+J(t,x,\mu_t,u)\}
\]
has the claimed properties.
\qed

To appreciate the results in Theorem \ref{thfeedbackHJB} better and to illustrate the conditions {\bf (A1)}, {\bf (A2)} and {\bf (A3)} are verifiable conditions for the feedback feedback
regularity property in \eqref{FBR}, we give the following example with concrete conditions.
\begin{example}
\label{example assumptions}
Set the control set $\U=[-1, 1]$. 
The controlled dynamics $X_t$ is associated to the family of operators $A$ of the form
\begin{align}
&A[t,\mu,u]f(x)\label{example A}\\
= &(h(t,x,\mu,u), \nabla f(x)) + L[t,\mu]f(x)\notag\\
=&(h(t,x,\mu,u), \nabla f(x)) + \frac{1}{2}(G(t,x,\mu)\nabla,\nabla) f(x)\notag
\end{align}
where the drift coefficient $h$ is linear in $u$ and of the form
$$h(t,x,\mu,u)=\int_{\R^d} \beta(t,x, y)\mu(dy) +u,$$
and 
$$G(t,x,\mu) = \int_{\R^d} g(t,x, y)\mu(dy)$$ 
with the functions $ \beta, g: [0,T]\times \R^d \times \R^d\to \R^d$ being bounded, continuous in $t$ and Lipschitz continuous uniformly in $x$ and $y$. See Remark \ref{example of h and L} for the corresponding form of stochastic differential equations.

Since the function $G$ is linear in $\mu$ hence differentiable in $\mu$. Together with the conditions on $g$,  the condition \eqref{eqassumonder2} is satisfied. 
The operator $L[t,\mu]f(x)=\frac{1}{2}(G(t,x,\mu)\nabla,\nabla) f(x)$ in \eqref{example A} generates a strongly continuous backward propagator which has the smoothing property with $\beta=\frac{1}{2}$, see discussion in subsection \ref{assumptions}. Hence, the assumption {\bf (A2)} is satisfied.

The running cost function $J$ is quadratic in $u$ and of the form
\begin{align}
&J(t,x,\mu,u)\notag\\
=&\int_{\R^d}\alpha (t,x,y)\mu(dy)-\frac{1}{2}u^2 \int_{\R^d}\theta(t,x,y)\mu(dy)\notag
\end{align}
where the functions $\alpha, \theta: [0,T]\times \R^d \times \R^d\to \R$ are bounded, continuous in $t$ and Lipschitz continuous uniformly in $x$ and $y$, and $\theta (t,x,y)> 0$ for any $(t,x,y)$. The boundeness of $J$  guarantees the condition \eqref{eq2thweak_existencem} is satisfied. The special form of the functions $h$ and $J$ guarantees the condition \eqref{eq1thweak_existencem} and \eqref{eqassumonder1} are satisfied. Hence  the assumption  {\bf (A1)} is satisfied. 

Together with  a  terminal function $V^T\equiv 0$,  Theorem \ref{Sensitivity 2} gives the result that the spacial derivative of the solution to \eqref{HJB_coupled} is Lipschitz continuous uniformly with respect to the functional parameter $\{\mu.\}$. 

Further, under this linear-quadratic setting, the max in \eqref{singleton control} is achieved only at one point and  one gets an explicit formula of the unique point of maximum, i.e.
{\small
 \begin{align}
&u(t,x,\mu,p)\label{formular u}\\
&=\arg \max_u \left\{ h(t,x,\mu,u) p+J(t,x,\mu,u)\right\}\notag\\
=&\begin{cases}
&\frac{p}{\int_{\R^d}\theta(t,x,y)\mu(dy)},\,\quad\quad\text{if }\frac{p}{\int_{\R^d}\theta(t,x,y)\mu(dy)}\in[-1,1]\notag\\[.7em]
  & 1, \quad\,\,\,\text{if } \int_{\R^d}\frac{\alpha}{\beta}(t,x,y)\mu(dy)p \notin[-1,1] \,\text{and } \,p>0\notag\\[.7em]
  & -1,\quad  \text{if } \int_{\R^d}\frac{\alpha}{\beta}(t,x,y)\mu(dy)p \notin[-1,1]\,\text{and } \,p<0\notag
  \end{cases}
\end{align}
}
which is continuous in $t$ and Lipschitz continuous in $(x,\mu, p)$ uniformly with respect to $x$, $\mu$ and bounded $p$.

Now it is checked that this example satisfies all the conditions in Theorem  \ref{thfeedbackHJB}. By setting $p=\nabla V(t,x;\{\mu.\})$, we may conclude that the result optimal feedback control strategy has the the feedback regularity property. 
\end{example}

\begin{remark}
Let us stress again that in \eqref{V Lip}, \eqref{V Lip_2} the space $\S$ is an abstract Banach space, but in application to control depending on empirical measures, we have in mind the norm of the dual space $(\C^2)^*$, where $\C^2$ is the domain of the generating family $A[t,\mu,u]$.
\end{remark}

%
%

\section{Appendix}
\appendix
\section{Proof of Theorem \ref{weak_existence} }
\label{Proof to the Theorem weak_existence }
Let $C_{V^T}^T([0,T], \C^1)$ be the set of functions $\phi: [0,T]\times \R^d\to \R$, which satisfy $\phi (T,x)=V^T(x)$ for all $x\in \R^d$, $\phi (t,\cdot)\in \C^1$ for each $t\in [0,T]$ with the norms $\|\phi(t,\cdot)\|_{\C^1}$  uniformly bounded in $t$, and the mapping $t\mapsto \phi(t,\cdot)$ is continuous as a mapping from $[0,T]$ to $\C$. We equip this space with the norm
$$\|\phi\|_{C_{V^T}^T([0,T], \C^1)}:=\sup_{t\in[0,T]}\|\phi(t,\cdot)\|_{\C^1}.$$ Note this definition of the set $C_{V^T}^T([0,T], \C^1)$ is not standard in the sense that it the continuity is considered from $[0,T]$ to $\C$, but not from $[0,T]$ to $\C^1$.

Define an operator $\Psi$ acting on $C_{V^T}^T([0,T], \C^1)$ by
\begin{align}
\label{mildformpsi1}
\Psi(\phi)(t,x):=& (U^{t,T} V^T(\cdot))(x) \notag\\
&+ \int_t^T U^{t,s}H_s(\cdot,\nabla \phi(s,\cdot))(x)ds.
\end{align}
Since the propagator $U^{t,T}$ is strongly continuous in $t$ and the integral term is continuous in $t$, the mapping $t\to \Psi(\phi)(t,\cdot)$ is continuous.

Since $V^T(\cdot)\in \C^1$ and the family $\{{U}^{t,T}, 0\leq t\leq T\}$ is bounded as a family of  mappings from $\C^1$ to $\C^1$, we have ${U}^{t,T}V^T(\cdot)\in \C^1$ and  it is uniformly bounded on $0\leq t\leq T$.
By the triangle inequality and \eqref{eq1thweak_existencem},\eqref{eq2thweak_existencem}, for each $t\in[0,T]$
{\small
\begin{eqnarray}
\label{ineq H}
&&\|H_t(\cdot,\nabla \phi(t,\cdot))\|_{\C}\notag\\
&\leq& \|H_t(\cdot,0)\|_{\C}
+\|H_t(\cdot,\nabla \phi(t,\cdot))-H_t(\cdot,0)\|_{\C} \nonumber \\[.4em]
              &\leq& c_2+c_1\|\nabla \phi(t,\cdot) \|_{\C} \nonumber \\[.4em]
              &\leq & c_2+c_1\| \phi(t,\cdot) \|_{\C^1}.
\end{eqnarray}
}
The last inequality in \eqref{ineq H} holds since by definition $\| \phi(t,\cdot) \|_{\C^1}= \|\nabla \phi(t,\cdot) \|_{\C}+ \|\phi(t,\cdot)\|_{\C}$ and $\|\phi(t,\cdot)\|_{\C}\geq 0$.
The smoothing condition \eqref{smoothingm} guarantees that
${U}^{t,s} H_s(\cdot,\nabla \phi(s,\cdot))  \in \C^1$ for each $0\le t<s\leq T$. The conditions \eqref{BDD 2m}and \eqref{smooth property 2m} and the inequality \eqref{ineq H} imply for each $t\in[0,T]$ that
{\small
\begin{align*}\label{GL}
&\|\Psi(\phi)(t,\cdot)\|_{\C^1}\notag\\
 \leq &\|{U}^{t,T}V^T(\cdot)\|_{{\C^1}}+\int_t^T \| U^{t,s} H_s(\cdot,\nabla \phi(s,\cdot))\|_{\C^1} \,\,ds \notag\\
\leq &c_4 \|V^T(\cdot) \|_{\C^1} + c_5 \int_t^T (s-t)^{-\be}\|H_s(\cdot,\nabla \phi(s,\cdot))\|_{\C}\,\, ds\notag\\
 \leq & c_4 \|V^T (\cdot)\|_{\C^1} + c_5 \int_t^T (s-t)^{-\be}(c_2+c_1\| \phi(s,\cdot) \|_{\C^1})\,\, ds\notag\\
  \leq & c_4 \|V^T (\cdot)\|_{\C^1} + c_5\left (c_2+c_1\sup_{t\leq s\leq T}\|\phi(s,\cdot)\|_{\C^1}\right)\frac{(T-t)^{1-\be}}{1-\be}.
\end{align*}
}

It follows that  the operator $\Psi$ maps $C_{V^T}^T([0,T], \C^1)$ to itself, i.e.
\[
\Psi: C_{V^T}^T([0,T], \C^1)\mapsto C_{V^T}^T([0,T], \C^1).
\]
Conditions \eqref{eq1thweak_existencem} and \eqref{smooth property 2m}  imply  for every $\phi^1,\phi^2 \in C_{V^T}^T([0,T],\C^1)$ and $t\in [0,T]$ that
{\small
\begin{align}
\label{eq5thweak_existence}
&\|\Psi(\phi^1)(t,\cdot)-\Psi(\phi^2)(t,\cdot)\|_{\C^1}\notag\\
\leq &\int_t^T \| U^{t,s}[H_s(\cdot,\nabla \phi^1(s,\cdot))-H_s(\cdot,\nabla \phi^2(s,\cdot))]\|_{\C^1}\,\,ds \notag \\
\leq &\int_t^T c_5 (s-t)^{-\be} c_1 \|\nabla \phi^1(s,\cdot)-\nabla \phi^2(s,\cdot) \|_{\C}\,\,ds \notag \\
\leq &c_1 c_5 \frac{(T-t)^{1-\be}}{1-\be} \sup_{t\leq s\leq T} \| \phi^1(s,\cdot)-\phi^2(s,\cdot) \|_{\C^1}.
\end{align}
}

Hence, for $T>0$ small enough, we get wellposedness by the contraction principle. Similar arguments yield the wellposedness on the interval $[t_0, T]$ for $T-t_0$ small enough, with $\Psi(\phi) (T-t_0, \cdot)\in \C^1$. Iteration of the above procedure on the whole interval $[0,T]$ completes the proof. 
\qed
\section{Proof of Theorem \ref{Sensitivity 1} }
\label{Proof to the Theorem Sensitivity 1}

From \eqref{smooth property 3}, for each $\al_1, \al_2\in[0,1]$, $t\in[0,T]$ and $x\in\R^d$, we have
\begin{align}\label{W-diff}
&W_{\al_1}(t,x)-W_{\al_2}(t,x)\\[.4em]
&=U_{\al_1}^{t,T}\left(  V^T_{\al_1} - V^T_{\al_2}\right)(x)+\left( U_{\al_1}^{t,T}-U_{\al_2}^{t,T}\right)V^T_{\al_2}(x).\notag
\end{align}

By the condition {\bf (A3)}, for each $x\in \R^d$ the mapping $\al\to V_{\al}^T(x)$ is differentiable and the derivative $\frac{\partial V^T_\al}{\partial \al}(\cdot)$ belongs to $\C^1$. Since for any $0\leq t\leq T$ and $\al_1\in[0,1]$, $U_{\al_1}^{t,T}: \C^1\to \C^1$, together with \eqref{BDD 2m} we have
\begin{align}\label{ex1}
&\|U_{\al_1}^{t,T}\left(  V^T_{\al_1} - V^T_{\al_2}\right)(\cdot)\|_{\C^1}\notag\\
&=\left\|U_{\al_1}^{t,T} \int_{\al_2}^{\al_1}\frac{\partial V_\gamma^T}{\partial \al}(\cdot)d\gamma\right\|_{\C^1}\notag\\
&\leq  c_4 |\al_1-\al_2|\sup_{\gamma\in[\al_1,\al_2]}\left\|\frac{\partial V_\gamma^T}{\partial \al}(\cdot)\right\|_{\C^1}.
\end{align}
By \eqref{PropergatorProperty } in Proposition \ref{prop-propergatorProperty} and the smoothing property \eqref{smooth property 2m}, we have
$$U_{\al_1}^{t,T}-U_{\al_2}^{t,T}=\int_t^T    U_{\al_1} ^{t,s}( L_{\al_1}[s]-L_{\al_2}[s])U_{\al_2} ^{s,T}ds$$
which is an operator mapping $\C^2$ to $\C^1$. Together with the condition {\bf (A2)} (ii) that for each $t\in[0,T]$ the mapping $\al\mapsto L_{\al}[t]$ is differentiable and $\frac{\partial L_\al}{\partial \al}[t]: \C^2\to\C$, we have
{\small
\begin{align}\label{ex2}
&\left\|\left( U_{\al_1}^{t,T}-U_{\al_2}^{t,T}\right)V^T_{\al_2}(\cdot)\right\|_{\C^1}\notag\\
&\leq c_4c_5 \int_t^T(s-t)^{-\beta}ds\notag\\
&\hspace{3em}\times\sup_{s\in[t,T]} \| L_{\al_1}[s]-L_{\al_2}[s]\|_{\C^2\to\C}\|V^T_{\al_2}(\cdot)\|_{\C^2}\notag\\
&\leq c_4c_5 \frac{(T-t)^{1-\be}}{1-\beta}|\al_1-\al_2|\notag\\
&\hspace{3em}\times\underset{\begin{subarray}{c}
  s\in[t,T] \\
  \gamma\in[\al_1,\al_2]
  \end{subarray}}{\sup}\left\| \frac{\partial L_\gamma}{\partial \al}[s]\right\|_{\C^2\to\C}\|V^T_{\al_2}(\cdot)\|_{\C^2}.
\end{align}
}
Therefore, from \eqref{W-diff} together with \eqref{ex1} and \eqref{ex2}, we complete the proof.
\qed

\section{Proof of Proposition \ref{prop}}
\label{Proof for the Proposition prop}
(i) Since the operator $L_\al[t]$ is differentiable in $\al$ for each $t\in[0,T]$, together with \eqref{PropergatorProperty } in Proposition 2.1, for $\al_1,\al\in[0,1]$ with $\al_1\neq\al$, we have
{\small
\begin{align*}\label{U-a}
\frac{U^{t,s}_{\al_1}-U^{t,s}_{\al}}{\al_1-\al}&=\frac{1}{\al_1-\al}\int^s_t U^{t,r}_{\al_1}(L_{\al_1}[r]- L_{\al}[r])U^{r,s}_{\al}dr\notag\\
&=\int^s_t U^{t,r}_{\al_1}\frac{\partial L_\gamma}{\partial \al}[r]U^{r,s}_{\al}dr
\end{align*}
}
with some $\gamma \in[\al_1,\al]$. By \eqref{PropergatorPropertya} in Proposition  \ref{prop-propergatorProperty},  for each $0\leq t<r\leq T$ the mapping $\al\mapsto U_{\al}^{t,r}$ is continuous as a function from $[0,1]$ to $\mathcal{L}(\C^2, \C)$ in the strong operator  topology.

Next, we approximate a continuous function by a sequence of twice continuously differentiable functions. 
Then by standard density arguments, we have for each $0\leq t<r\leq T$, the mapping $\al\mapsto U_{\al}^{t,r}$ is also continuous as a function from $[0,1]$ to $\mathcal{L}(\C, \C)$ in the strong operator  topology.

Together with the smoothing property \eqref{smoothingm} and the condition that for each $r\in[0,T]$ and $\gamma\in[0,1]$, $\frac{\partial L_\gamma}{\partial \al}[r]: \C^2\to \C$, we have that the derivative
$$
\frac{\partial U^{t,r}_\al}{\partial \al}:=\lim_{\al_1\to\al}\frac{U^{t,s}_{\al_1}-U^{t,s}_{\al}}{\al_1-\al}
$$
 exists in the strong topology in $\mathcal{L}(\C^2, \C)$ and equals to
$$\frac{\partial U^{t,r}_\al}{\partial \al}=\int^s_t U^{t,r}_{\al}\frac{\partial L_\al}{\partial \al}[r]U^{r,s}_{\al}dr.$$

(ii) By condition {\bf (A3)}, the mapping $\al\mapsto V^T_\al(\cdot)$ is differentiable and the derivative exists in $\C^1$, namely, $\lim_{\al_1\to\al}\frac{V^T_{\al_1}-V^T_{\al}}{\al_1-\al}$ exists and belongs to $\C^1$. Since the mapping $\al\mapsto U_{\al}^{t,r}$ is continuous as a function from $[0,1]$ to $\mathcal{L}(\C, \C)$ in the strong operator  topology,  hence for any $t\in[0,T]$,
$$\lim_{\al_1\to\al}\left(U^{t,T}_{\al_1}\frac{V^T_{\al_1}-V^T_{\al}}{\al_1-\al}\right)=U^{t,T}_{\al}\frac{V^T_{\al}}{\partial \al}\in\C.$$
 Finally,  from \eqref{smooth property 3} , the derivative
\begin{align}
\frac{\partial W_{\al}}{\partial \al}(t,\cdot)&:=\lim_{\al_1\to\al}\frac{W_{\al_1}(t,\cdot)-W_{\al}(t,\cdot)}{\al_1-\al}\notag\\
&=\lim_{\al_1\to\al}\left(U^{t,T}_{\al_1}\frac{V^T_{\al_1}(\cdot)-V^T_{\al}(\cdot)}{\al_1-\al}\right)\notag\\
&\hspace{1em}+ \lim_{\al_1\to\al}\frac{U^{t,T}_{\al_1}-U^{t,T}_{\al}}{\al_1-\al}V^T_{\al}(\cdot)\notag\\
&={U}_\alpha^{t,T} \frac{\partial V^T_\alpha}{\partial \alpha} (\cdot)
+\frac{\partial U_\alpha^{t,T}}{\partial \alpha}V^T_\alpha(\cdot)\notag
\end{align}
exists in $\C$ for each $t\in[0,T]$ and any $\al\in[0,1]$.
\qed

\section{Proof of Theorem \ref{Sensitivity 2}}
\label{Proof for the Theorem Sensitivity 2}

(a) Recall in the proof of Theorem \ref{weak_existence}, for any $\alpha \in [0,1]$,  the unique solution $V_\alpha $ is the unique fixed point of the mapping
$$\phi \mapsto \Psi_{\al}(\phi), \quad C^T_{V_\al^T}([0,T], \C^1)\to C^T_{V_\al^T}([0,T], \C^1)$$
defined by, for each $t\in[0,T]$
\begin{equation}
\label{mapping Psi}
\Psi_{\al}(\phi)(t,\cdot) = U_\alpha^{t,T} V_\alpha^T(\cdot) + \int_t^T U_\alpha^{t,s}H_{\alpha,s}(\cdot,\nabla \phi (s,\cdot))ds.
\end{equation}
For any $\al_i\in[0,1]$, $i=1,2$, let $V_{\al_i}$ be the unique fixed point of the mapping $\Psi_{\al_i}$, i.e.
$$V_{\al_i}=\Psi_{\al_i} (V_{\alpha_i} ),\quad \text{for } i=1,2.$$
Then from \eqref{mapping Psi} we have
{\small
\begin{align}
&V_{\al_1}(t, \cdot)-V_{\al_2}(t, \cdot)\notag\\
=&\Phi_{\al_1}(V_{ \al_1}(t, \cdot))-\Phi_{\al_2}(V_{\al_2}(t, \cdot))\notag\\[.4em]
 =&U^{t,T}_{\al_1}V^T_{\al_1}(\cdot)-U^{t,T}_{\al_2}V^T_{\al_2}(\cdot)\notag\\[.4em]
+&\int_t^TU^{t,s}_{\al_1}H_{\al_1,s}(\cdot,\nabla V_{\al_1}(s,\cdot))ds \notag\\
-&\int_t^TU^{t,s}_{\al_2}H_{\al_2,s}(\cdot,\nabla V_{\al_2}(s,\cdot))ds\notag\\[.4em]
=&U^{t,T}_{\al_1}V^T_{\al_1}(\cdot)-U^{t,T}_{\al_2}V^T_{\al_2}(\cdot)\notag\\[.4em]
+&\int_t^T \left(U^{t,s}_{\al_1}- U^{t,s}_{\al_2}\right)H_{\al_1,s}(\cdot,\nabla V_{\al_1}(s,\cdot))ds \notag\\
+&\int_t^T U^{t,s}_{\al_2}\left(H_{\al_1,s}(\cdot,\nabla V_{\al_1}(s,\cdot))-H_{\al_2,s}(\cdot,\nabla V_{\al_1}(s,\cdot))\right)ds \notag\\
+&\int_t^T U^{t,s}_{\al_2}\left(H_{\al_2,s}(\cdot,\nabla V_{\al_1}(s,\cdot))-H_{\al_2,s}(\cdot,\nabla V_{\al_2}(s,\cdot))\right)ds\notag\\
=:&\Lambda_1+\Lambda_2+\Lambda_3+\Lambda_4 \label{5}.
\end{align}
}
By theorem \ref{Sensitivity 1}, there exists a constant $c>0$ such that
{\small
\begin{align}\label{Lambda_1}
&\|\Lambda_1\|_{\C^1}
\leq c |\al_1-\al_2|  \sup_{\gamma\in[\al_1,\al_2]}\left\|\frac{\partial V^T_{\gamma}}{\partial \al}\right\|_{\C^1}\notag\\
&+c |\al_1-\al_2|
(T-t)^{1-\beta}\hspace{-.8em}\underset{\begin{subarray}{c}
  s\in[t,T] \\
  \gamma\in[\al_1,\al_2]
  \end{subarray}}{\sup}\left\| \frac{\partial L_\gamma}{\partial \al}[s]\right\|_{\C^2\to\C}\|V^T_{\al_2}(\cdot)\|_{\C^2}.
\end{align}
}
By proposition \ref{prop-propergatorProperty} and the differentiability of $L_\al[t]$ in $\al$ for each $t\in[0,T]$, we have
{\small
\begin{align}\label{Lambda_2}
&\|\Lambda_2\|_{\C^1}\\
&=\Big\|\int_t^T \left(     \int_t^s U^{t,r}_{\al_1}(L_{\al_1}[r]-L_{\al_2}[r])U^{r,s}_{\al_2}dr\right)\notag\\
&\hspace{3cm}H_{\al_1,s}(\cdot,\nabla V_{\al_1}(s,\cdot))ds\Big\|_{\C^1} \notag\\[.4em]
&=\Big\| \int_t^T \left(     \int_t^s U^{t,r}_{\al_1}(\al_1-\al_2)\frac{\partial L_\theta}{\partial \al}[r]U^{r,s}_{\al_2}dr\right) \notag\\[.4em]
&\hspace{3cm}H_{\al_1,s}(\cdot,\nabla V_{\al_1}(s,\cdot))ds\Big\|_{\C^1} \notag\\[.4em]
& \leq|\al_1-\al_2|c_5c_6\int_t^T\int_t^s(r-t)^{-\beta}(s-r)^{-\beta}\|  \frac{\partial L_\theta}{\partial \al}[r] \|_{\C^2\to\C}dr\notag\\
&\hspace{3cm} \| H_{\al_1,s}(\cdot,\nabla V_{\al_1}(s,\cdot)) \|_{\C_{Lip}}ds\notag\\[.4em]
 & \leq|\al_1-\al_2|c_5c_6 m\frac{ (T-t)^{2-2\be}}{(1-\be)^2} \sup_{r\in[t,T]}\left\|\frac{\partial L_\theta}{\partial \al}[r]\right\|_{\C^2\to \C}\notag
\end{align}
}
with $\theta\in[0,1]$ and $m:=\sup_{\alpha\in[0,1]}\sup_{( s,p)\in \O}\|H_{\al,s}(\cdot, p)\|_{\C_{Lip}}<\infty$. The last inequality in \eqref{Lambda_2} is obtained through the calculation
\begin{align*}
&\int_t^T \int_t^s (r-t)^{-\beta}(s-r)^{-\beta}drds\\
=&\int_t^T \int_r^T (r-t)^{-\beta}(s-r)^{-\beta}dsdr\\
=&  \int_t^T (r-t)^{-\beta} \frac{(T-r)^{1-\beta}}{1-\beta}dr\\
\leq &\frac{(T-t)^{1-\beta}}{1-\beta}\int _t^T(r-t)^{-\beta}dr=\frac{ (T-t)^{2-2\be}}{(1-\be)^2}.
\end{align*}

By the condition {\bf (A1)}, for each $t\in[0,T]$ and $x,p\in\R^d$ the mapping $\al\mapsto H_{\al,s}(x,p)$ is differentiable and the derivative is continuous, we have
{\small
\begin{align}\label{Lambda_3}
&\|\Lambda_3\|_{\C^1}\notag\\
&=\left\|\int_t^T U^{t,s}_{\al_2} (\al_1-\al_2) \frac{\partial H_{\al,s}}{\partial \al}(\cdot,\nabla V_{\al_1}(s,\cdot))ds\right \|_{\C^1}\notag\\
&\leq |\al_1-\al_2|\int_t^T c_5 (s-t)^{-\be}\left\|\frac{\partial H_{\theta,s}}{\partial \al}(\cdot, \nabla V_{\al_1}(s,\cdot))\right\|_{\C} ds\notag\\[.4em]
&\leq  |\al_1-\al_2| c_5 \frac{(T-t)^{1-\be}}{1-\be} \,  \sup_{(s,p)\in\O}\left \|\frac{\partial H_{\theta,s}}{\partial \al}(\cdot, p)\right\|_{\C}.
\end{align}
}
By \eqref{smooth property 2m} and \eqref{eq1thweak_existencem}, we get
{\small
\begin{align}\label{Lambda_4}
&\|\Lambda_4\|_{\C^1}\notag\\
=&\left\| \int_t^T \hspace{-.6em}U^{t,s}_{\al_2}\left(H_{\al_2,s}(\cdot,\nabla V_{\al_1}(s,\cdot))-H_{\al_2,s}(\cdot,\nabla V_{\al_2}(s,\cdot))\right)ds \right\|_{\C^1}\notag\\
\leq &c_1 c_5\int_t^T (s-t)^{-\beta} \|\nabla V_{\al_1}(s,\cdot)-\nabla V_{\al_2}(s,\cdot)\|_{\C}  \,        ds\notag\\
\leq &c_1 c_5\int_t^T (s-t)^{-\beta} \|V_{\al_1}(s,\cdot)- V_{\al_2}(s,\cdot)\|_{\C^1}  \,      ds\notag\\
\leq &c_1 c_5\frac{(T-t)^{1-\beta}}{1-\beta} \sup_{s\in[t,T]}\|V_{\al_1}(s,\cdot)- V_{\al_2}(s,\cdot)\|_{\C^1}.
\end{align}
}

It follows, from (\ref{5}) together with the estimates \eqref{Lambda_1}, \eqref{Lambda_2},  \eqref{Lambda_3} and  \eqref{Lambda_4}, that
{\small
\begin{align}
&\sup_{t\in[0,T]}\| V_{\al_1}(t,\cdot)-V_{\al_2}(t,\cdot)\|_{\C^1}\\
&\leq c |\al_1-\al_2|  \sup_{\gamma\in[\al_1,\al_2]}\left\|\frac{\partial V^T_{\gamma}}{\partial \al}\right\|_{\C^1}\notag\\
&+c |\al_1-\al_2|
(T-t)^{1-\beta}\|V^T_{\al_2}(\cdot)\|_{\C^2}
\underset{\begin{subarray}{c}
  t\in[0,T] \\
  \gamma\in[\al_1,\al_2]
  \end{subarray}}{\sup}\left\| \frac{\partial L_\gamma}{\partial \al}[t]\right\|_{\C^2\to\C}\notag\\
&+|\al_1-\al_2|c_5c_6 m\frac{(T-t)^{2-2\be} }{(1-\be)^2} \underset{\begin{subarray}{c}
  t\in[0,T] \\
  \gamma\in[\al_1,\al_2]
  \end{subarray}}{\sup}\left\|\frac{\partial L_\gamma}{\partial \al}[t]\right\|_{\C^2\to \C}\notag\\
&+ |\al_1-\al_2| c_5 \frac{(T-t)^{1-\be}}{1-\be} \,
\underset{\begin{subarray}{c}
  (t,p)\in\O \\
  \gamma\in[\al_1,\al_2]
  \end{subarray}}{\sup}
\left \|\frac{\partial H_{\gamma,t}}{\partial \al}(\cdot, p)\right\|_{\C}\notag\\
&+c_1 c_5\frac{(T-t)^{1-\beta} }{1-\beta}\sup_{t\in[0,T]}\|V_{\al_1}(t,\cdot)- V_{\al_2}(t,\cdot)\|_{\C^1} .
\end{align}
}
For $t$ close to $T$ enough so that $\frac{(T-t)^{1-\beta}}{(1-\beta)^2} <1$, we have inequality \eqref{T10}. For arbitrary $t\in[0,T]$, the proof follows by iterations.

(b)
By the definitions of $L_\al[t]$, $H_{\al,t}(x,p)$, $V^T_\al(x)$ in \eqref{Halpha}, \eqref{Lalpha}, \eqref{VTalpha} respectively and the assumptions \eqref{eqassumonder1}, \eqref{eqassumonder2}, \eqref{eqassumonder3}, for any $\{\mu^1_.\},\{\mu^2_.\}\in C([0,T], \M)$,
the statement follows from the equation \eqref{estimate} and the inequality \eqref{T10} by setting $\al_1=1$ and $\al_2=0$.
\qed

{\bf Acknowledgements: }
The authors thank the referees for their valuable suggestions.

\end{document}